\documentclass[12pt]{amsart}
\usepackage{amssymb,amsmath}
\numberwithin{equation}{section}
\def\T{\text}
\def\simgeq{\underset\sim>}
\def\simleq{\underset\sim<}
\def\1#1{\overline{#1}}
\def\2#1{\widetilde{#1}}
\def\3#1{\widehat{#1}}
\def\4#1{\mathbb{#1}}

\def\5#1{\frak{#1}}
\def\6#1{{\mathcal{#1}}}
\def\C{{\4C}}
\def\R{{\4R}}

\def\Z{{\4Z}}

\begin{document}
\abstract
We prove that CR lines in an exponentially degenerate boundary are propagators of holomorphic extension. This explains, in the context of the CR geometry, why in this situation the induced Kohn-Laplacian $\Box_b$   is not hypoelliptic (Christ \cite{Ch02}).
\endabstract
\title[Propagation of regularity]{Propagation of regularity for solutions of the Kohn Laplacian in a flat boundary}
\author[L.~Baracco, T.V.~Khanh and G.~Zampieri]{Luca Baracco, Tran Vu Khanh and  Giuseppe Zampieri}
\address{Dipartimento di Matematica, Universit\`a di Padova, via 
Trieste 63, 35121 Padova, Italy}
\email{baracco@math.unipd.it, khanh@math.unipd.it,  
zampieri@math.unipd.it}
%\subjclass{}
\maketitle
%\tableofcontents
% Standard sets
\def\Giialpha{\mathcal G^{i,i\alpha}}
\def\cn{{\C^n}}
\def\cnn{{\C^{n'}}}
\def\ocn{\2{\C^n}}
\def\ocnn{\2{\C^{n'}}}
% Abbreviations
\def\const{{\rm const}}
\def\rk{{\rm rank\,}}
\def\id{{\sf id}}
\def\aut{{\sf aut}}
\def\Aut{{\sf Aut}}
\def\CR{{\rm CR}}
\def\GL{{\sf GL}}
\def\Re{{\sf Re}\,}
\def\Im{{\sf Im}\,}
\def\codim{{\rm codim}}
\def\crd{\dim_{{\rm CR}}}
\def\crc{{\rm codim_{CR}}}
\def\phi{\varphi}
\def\eps{\varepsilon}
\def\d{\partial}
\def\a{\alpha}
\def\b{\beta}
\def\g{\gamma}
\def\G{\Gamma}
\def\D{\Delta}
\def\Om{\Omega}
\def\k{\kappa}
\def\l{\lambda}
\def\L{\Lambda}
\def\z{{\bar z}}
\def\w{{\bar w}}
\def\Z{{\1Z}}
\def\t{{\tau}}
\def\th{\theta}
\emergencystretch15pt
\frenchspacing
\newtheorem{Thm}{Theorem}[section]
\newtheorem{Cor}[Thm]{Corollary}
\newtheorem{Pro}[Thm]{Proposition}
\newtheorem{Lem}[Thm]{Lemma}
\theoremstyle{definition}\newtheorem{Def}[Thm]{Definition}
\theoremstyle{remark}
\newtheorem{Rem}[Thm]{Remark}
\newtheorem{Exa}[Thm]{Example}
\newtheorem{Exs}[Thm]{Examples}
\def\Label#1{\label{#1}}
\def\bl{\begin{Lem}}
\def\el{\end{Lem}}
\def\bp{\begin{Pro}}
\def\ep{\end{Pro}}
\def\bt{\begin{Thm}}
\def\et{\end{Thm}}
\def\bc{\begin{Cor}}
\def\ec{\end{Cor}}
\def\bd{\begin{Def}}
\def\ed{\end{Def}}
\def\br{\begin{Rem}}
\def\er{\end{Rem}}
\def\be{\begin{Exa}}
\def\ee{\end{Exa}}
\def\bpf{\begin{proof}}
\def\epf{\end{proof}}
\def\ben{\begin{enumerate}}
\def\een{\end{enumerate}}
\def\dotgamma{\Gamma}
\def\dothatgamma{ {\hat\Gamma}}

\def\simto{\overset\sim\to\to}
\def\1alpha{[\frac1\alpha]}
\def\T{\text}
\def\R{{\Bbb R}}
\def\I{{\Bbb I}}
\def\C{{\Bbb C}}
\def\Z{{\Bbb Z}}
\def\Fialpha{{\mathcal F^{i,\alpha}}}
\def\Fiialpha{{\mathcal F^{i,i\alpha}}}
\def\Figamma{{\mathcal F^{i,\gamma}}}
\def\Real{\Re}
%
%
%
%\endtopmatter
\section{  Introduction}
J.J. Kohn noticed in \cite{K77} that analytic discs in the boundary of a pseudoconvex domain $\Omega\subset\C^n$ prevent from the $C^\infty$-hypoellipticity of the $\bar\partial$-Neumann problem: the canonical solution is not smooth exactly at the boundary points where the datum is. On the other hand, it has been explained by N.~Hanges and F.~Treves in \cite{HT83} that discs seating in $\partial\Omega$ are propagators of holomorphic extension from $\Omega$ across $\partial\Omega$. Thus, propagation and hypoellipticity appear to be in contrast one to another. M.~Christ proved in \cite{Ch02} that in the hypersurface in $\C^2$ defined by
\begin{equation}
\Label{1.1}
x_2=e^{-\frac1{|y_1|^s}},
\end{equation}
 one does not have hypoellipticity for the  induced Kohn-Laplacian $\Box_b$ when $s\geq1$. Note that for $s<1$ 
 this
  is hypoelliptic as well as the $\bar\partial$-Neumann problem: in fact, in this case, one has superlogarithmic estimates which are sufficient for hypoellipticity. It is worth remarking that superlogarithmicity does not entirily rule hypoellipticity. The pseudoconvex domain whose boundary 
is defined by the same equation as \eqref{1.1} but with $y_1$ replaced by $z_1$, that is, $x_2= e^{-\frac1{|z_1|^s}}$, has the same range $s<1$ for superlogarithmic estimates and, nonetheless, there always is hypoellipticity, for any value of $s$ (Kohn \cite{K00}). Here the matter is of a genuinily geometric type: there are no  curves running in complex tangential directions along which the manifold is flat and which are, therefore, possible propagators. Coming back to the pseudoconvex domain with boundary \eqref{1.1}, we show here that the lines parallel to the $x_1$-axis are propagators of holomorphic extension when $s\geq1$.
 More precisely, our statement is that discs in $\C^2$ over the 1-dimensional discs squeezed along the these lines, singular at $x_1=0$ and with boundary in $\partial\Omega$ apart from $x_1=1$ where they enter in $x_2<0$, ``point down" at $x_1=0$ if and only if $s\geq1$. These discs propagate the extendibility down; 
 in particular, this cannot be proved when $s<1$.
There is no surprise about it because, for $s<1$, these lines are not propagators of smoothness at the boundary. This follows from the hypoellipticity of the $\bar\partial$-Neumann problem. In fact, let $\chi=\chi(x_1)$ be $C^\infty$ and   satisfy $\chi\equiv0$ at $1$ and $\chi\equiv 1$ at $1$, and consider the $\bar\partial$-closed form $f:=\bar\partial\left(\frac{\chi(x_1)}{z_2}\right)$. If $s<1$ the equation $\bar\partial u=f$ has a solution $u$ in $\Omega$ which is smooth at $0$ and $1$; 
 thus the difference $u-\frac{\chi(x_1)}{z_2}$ is holomorphic in $\Omega$, singular at $x_1=0$ but smooth at $x_1=1$. 

We are indebted to Alexander Tumanov for important advice.

\section{Squeezing discs along lines}
In the standard disc $\Delta$ of the complex plane $\C$ with variable $\tau=re^{i\theta},\,\,\theta\in[0,2\pi]$, we consider the family of holomorphic mappings (=discs) depending on a small real parameter $\alpha$:
$$
\phi_\alpha(\tau)=-\frac1{\log(\frac14\left(\frac{1-\tau}2\right)^\alpha)}.
$$
The discs are squeezed along the interval $(0,|\log\frac14|^{-1})$ as $\alpha\searrow0$ with the points $+1$ and $-1$ interchanged with the left and right bounds respectively and they are singular at $\tau=1$. Moreover, the most of their mass concentrates at $\tau=-1$. 
We have
\begin{equation*}
\frac1{|\phi_\alpha(\tau)|}\sim -\alpha\log\left(\frac{|1-\tau|}2\right),\quad\tau\in\Delta.
\end{equation*}
With the notation $\tau=e^{i\theta}\in\partial\Delta$ we also have
\begin{equation*}
\begin{split}
\arg\left(\frac{1-\tau}2\right)&=\T{arctg}\left(\frac{\sin\theta}{1-\cos \theta}\right)
\\
&=\T{arctg}\left(\frac{\cos\frac\theta2\sin\frac\theta2}{\sin^2\frac\theta2}\right)
\\
&=\T{arctg}\left(\T{cotg}\frac\theta2\right)=\frac\pi2-\frac\theta2
\end{split}
\end{equation*}
and finally
\begin{equation*}
\begin{split}
\frac1{|\Im\phi_\alpha|}&\sim\frac{\log^2\left(\frac14\frac{|1-\tau|^\alpha}2\right)+\alpha^2(\frac\pi2-\frac\theta2)^2}{\alpha|\frac\pi2-\frac\theta2|}
\\
&\sim\frac{\log^2(\frac14\frac{|1-\tau|^\alpha}2)}\alpha+O(\alpha)
\\
&=\frac1\alpha\log^2\frac14+2\log\frac14\log\frac{|1-\tau|}2+\alpha\log^2\frac{|1-\tau|}2+O(\alpha).
\end{split}
\end{equation*}
Thus, for $\alpha$ fixed as in next proposition, we have $\frac1{|\Im\phi_\alpha|}\sim\alpha\log^2\frac{|1-\tau|}2$ at $\tau=1$.
\bp
\Label{p1.1}
We have 
\begin{itemize}
\item[(i)] $ e^{-\frac1{|\phi_\alpha|^s}}=O^\infty(1-\tau)$  for $s>1$,
\item[(ii)] $e^{-\frac1{|\Im\phi_\alpha|^s}}=O^\infty(1-\tau)$ for $s>\frac12$.
\end{itemize}
\ep
\bpf
As for (i), we have to notice that 
\begin{equation*}
\begin{split}
e^{-\frac1{|\phi_\alpha|^s}}&\sim e^{-\alpha\log^s\left(\frac{|1-\tau|}2\right)}
\\
&=|1-\tau|^{\alpha\left|\log^{s-1}\left(\frac{|1-\tau|}2\right)\right|}=O^\infty(|1-\tau|)\quad\T{for $s>1$}.
\end{split}
\end{equation*}
As for (ii), this follows from
\begin{equation*}
\begin{split}
e^{-\frac1{|\Im\phi_\alpha|^s}}&\sim e^{\alpha\log^{2s}\left(\frac{|1-\tau|}2\right)}
\\
&=|1-\tau|^{\alpha\log^{2s-1}\left|\frac{1-\tau}2\right|}=O^\infty(|1-\tau|)\quad\T{for $s>\frac12$}.
\end{split}
\end{equation*}
This concludes the proof of the proposition.

\epf
In particular, the two functions in the statement of the proposition are $C^\infty$, and thus also $C^{1,\beta}$, at $\tau=1$ for $s>1$ and $s>\frac12$ in the two  respective cases. We have a basic result about composition of $\phi_\alpha$ with flat functions more general than $e^{-\frac1{|z_1|^s}}$ or $e^{-\frac1{|y_1|^s}}$. For this, let $h_\eta(z_1,y_2),\,(z_1,y_2)\in\C\times\R$, be a function sufficiently smooth depending on a parameter $\eta$. 
\bp
\Label{p1.2}
Let $\eta\mapsto h_\eta,\,\,\R\to\C^3$ be $C^k$ and satisfy $\partial_\eta h_\eta\equiv0$ in a neighborhood of $z_1=0$. Assume that all (mixed) derivatives up to order 2 in $\tau$ and $k$ in $\eta$ are $O(e^{-\frac1{|y_1|^s}})$ for $s\geq\frac12$. Then, the function $(\eta,v)\mapsto h_\eta(\phi_\alpha,v)$ has the properties:
\begin{itemize}
\item[(i)] it sends $\R\times C^{1,\beta}\to C^{1,\beta}$,
\item[(ii)] it is $C^k$ with respect to $\eta$,
\item[(iii)] it is differentiable with respect to $v$ at $v=0$ and its differential is close to $0$. 
\end{itemize}
\ep
\bpf
(i): For a function $g$ of a real variable $t$, the assumptions
$$
g=O(e^{-\frac1{t^s}}),\quad g'=O(e^{-\frac1{t^s}}),
$$
imply
\begin{itemize}
\item[] $g(|\Im \phi_\alpha|)=O^\infty(t)$ when $s>\frac12$ (Proposition~\ref{p1.1}),
\item[] $|\Im\phi_\alpha|'\leq \frac1{|\log^3(1-\tau)|}\frac1{|1-\tau|}$,
\item[] $(g(|\Im \phi_\alpha|))'=g'(\Im\phi_\alpha)\frac1{|\log^3(1-\tau)||1-\tau|}=0^\infty(|1-\tau|)$ (again, Proposition~\ref{p1.1}).
\end{itemize}
This concludes the proof of (i).

(ii): Since $\partial_\eta h_\eta\equiv0$ when $\phi_\alpha$ is singular, then the $C^k$ dependence of $h_\eta(\phi_\alpha,v)$ on $\eta$ is a standard fact: if $\eta\mapsto g_\eta,\quad \R\to C^2(\R)$ is $C^2$ and $\sigma\in C^{1,\beta}$, then $\eta\mapsto g_\eta(\sigma),\quad \R\to C^{1,\beta}$ is $C^k$.

(iii): It is convenient to use a more general setting. Thus, let $g_\eta$ be $C^3$. Then, $v\overset {G_\eta}\mapsto g_\eta(v),\quad C^{1,\beta}\to C^{1,\beta}$ is $C^1$ at $v=0$ and its differential satisfies
$$
||G'_\eta|_{v=0}||_{L(C^{1,\beta},C^{1,\beta})}\simleq ||g_\eta||_{C^3}.
$$
Note that, in our application, $g_\eta=h_\eta(\phi_\alpha,\cdot)$; thus $||g_\eta||_{C^3}$ is small near $v=0$.

\epf
Now, we can set up a Bishop's equation in the unknown $v\in C^{1,\beta}$
\begin{equation}
\Label{bishop}
v-T_1(h_\eta(\phi_\alpha,v)=0,
\end{equation}
where $T_1$ is the Hilbert transform normalized by taking value $0$ at $\tau=1$. We rewrite the equation \eqref{bishop} in the functional space $C^{1,\beta}$ as $G_\eta(v)=0$. By (iii) of Proposition~\ref{p1.2}, we have
\begin{equation}
\Label{differential}
||G'_\eta
\big|_{v=0}-\T{id}||_{L(C^{1,\beta},C^{1,\beta})}\simleq ||h_\eta(\Im\phi_\alpha,0)||_{C^3}.
\end{equation}
By the implicit function theorem, we readily get
\bc
\Label{c1.1}
For small $\eta$, the equation \eqref{bishop} has a unique solution $v\in C^{1,\beta}$ and this depends in a $C^k$-fashion on $\eta$.
\ec
We write $v=v_{\alpha,\eta}$ for the solution of \eqref{bishop} and also write $u=-T_1v$ and $u=u_{\alpha,\eta}$. We also denote by $A=A_{\alpha\,\eta}$ the disc $A=(\phi,u+iv)$. When only dependence on $\alpha$ is relevant, we write $v=v_\alpha$, $u=u_\alpha$ and $A=A_\alpha$.  Note that under our assumption $h=O(e^{-\frac1{|y_1|^s}})$  we have 
$$
u=O(e^{-\frac1{|y_1|^s}}).
$$
For $\tau=re^{i\theta}$ and for a function in $C^0(\partial\Delta)$, such as $u_\alpha$, the harmonic extension of $u_\alpha$ from $\partial\Delta$ to $\Delta$, that we still denote by $u_\alpha$, has a radial derivative which is given by
\begin{equation}
\Label{radial}
\partial_ru_\alpha|_{\tau=1}=-p.v.\int_0^{2\pi}\frac{u_\alpha}{1-\cos\theta}d\theta,
\end{equation}
where the integral is taken in the sense of the principal value.
We first show that the values of $\theta$ for which $\phi_\alpha$ is not contained in the $\delta$-neighborhood of $|\log\frac14|^{-1}$ is very small.
\bl
\Label{l1.1}
We have the iclusion
\begin{equation}
\Label{**}
\left\{\theta:\,\left|\phi_\alpha(\theta)-|\log\frac14|^{-1}\right|>\delta\right\}\subset[0,e^{-\frac\delta\alpha}]\cup[2\pi-e^{-\frac\delta\alpha},2\pi],
\end{equation}
that is, the two intervals in the right side of \eqref{**} are sent, via $\phi_\alpha$, into the $\delta$-neighborhood of $\phi_\alpha(-1)$.
\el
\bpf
We have
$$
\left|-\frac1{\log(\frac14\left(\frac{1-\tau}2\right|^\alpha)}+\frac1{\log\frac14}\right|=\left|\frac{\log\left(\frac{1-\tau}2\right)^\alpha}{\log^2\frac14+\log\frac14\log\left(\frac{1-\tau}2\right)^\alpha}\right|.
$$
Now, the denominator is bounded away from $0$. Hence, the set in the left of \eqref{**} is contained in
$$
\{\tau:\,\left|\log\left(\frac{1-\tau}2\right)^\alpha\right|<\frac\delta\alpha\},
$$
which is in turn contained in the set
$$
\{\tau=e^{i\theta}:\,\theta<e^{-\frac\delta\alpha}\quad\T{or} \quad \theta>2\pi-e^{-\frac\delta\alpha}\}.
$$

\epf
Taking into account of Lemma~\ref{l1.1}, 
we decompose the integration in \eqref{radial} as 
$$
\partial_ru_{\alpha}=-\int_0^{2\pi}\cdot=-2\int_0^{e^{-\frac\delta\alpha}}\cdot-2\int_{e^{-\frac\delta\alpha}}^\pi\cdot.
$$
We approximate, near $\theta=0$, $1-\cos\theta$ by $\theta^2$ and define
$$
F_\alpha:=\int_0^{e^{-\frac\delta\alpha}}\frac{e^{-\frac1{|\Im\phi_\alpha|^s}}}{\theta^2}d\theta.
$$
\bp
\Label{p1.3}
(i) For $s\geq1$, we have $\underset{\alpha\to0}\lim F_\alpha=0$.

(ii) For $s<1$, we have $\underset{\alpha\to0}\lim F_\alpha=+\infty$.
\ep
\bpf
For the purpose of this proof it is not restrictive to replace $\frac\delta\alpha$ by $\frac1\alpha$. 

(i): Note that $\left|\frac{1-\tau}2\right|\sim\theta$ on the unit circle near $\tau=1$. We have
\begin{equation*}
\begin{split}
F_\alpha&\leq \int_0^{e^{-\frac1\alpha}}\frac{e^{-\frac1{|\Im\phi_\alpha|}}}{\theta^2}d\theta\quad\T{(since $s\geq1$})
\\
&\sim\int_0^{e^{-\frac1\alpha}}\frac{e^{-[\frac1\alpha\log^2\frac14+2\log\frac14\log\theta+\alpha\log^2\theta]}}{\theta^2}d\theta
\\
&\leq\int_0^{e^{-\frac1\alpha}}\frac{e^{-2\log\frac14\log\theta}}{\theta^2}d\theta
\\
&\leq\int_0^{e^{-\frac1\alpha}}1\,d\theta\leq {e^{-\frac1\alpha}}.
\end{split}
\end{equation*}
This proves (i).

(ii): We assume now $s<1$ and also suppose, without loss of generality, $s>\frac12$. By using the substitution $-\log\theta=t$, we get
\begin{equation*}
\begin{split}
F_\alpha&\geq \int_0^{e^{-\frac1\alpha}}e^{-\alpha^s\log^{2s}\theta-2\log\theta}d\theta
\\
&=\int_{\frac1\alpha}^{+\infty}e^{-\alpha^st^{2s}+2t}dt.
\end{split}
\end{equation*}
Now, we remark that $-\alpha^st^{2s}+2t>0$ if and only if $t<(\frac{2^{\frac1s}}\alpha)^{\frac s{2s-1}}$. Thus,
\begin{equation*}
\begin{split}
\int_0^{+\infty}\cdot&\geq \int_{\frac1\alpha}^{(\frac1\alpha)^{\frac s{2s-1}}}1dt
\\
&\simgeq {\left(\frac{2^{\frac1s}}\alpha\right)^{\frac s{2s-1}}}-\frac1\alpha\to+\infty,
\end{split}
\end{equation*}
where the last conclusion follows from $\frac s{2s-1}>1$.

\epf
\bt
\label{t1.1}
Let $\Omega\subset\C^2$ be a domain defined by $x_2>h(z_1,y_2)$ with $h$ satisfying $\partial^j_{z_1}h=O(e^{-\frac1{|y_1|^s}})$ with $s\geq1$ for any $j\leq 2$. Assume $s\geq 1$. Then, the lines $L\subset \partial\Omega$ defined by $y_1=0,\,\,x_2=\T{const}$ are propagators of holomorphic extendibility. Namely, if $f\in\T{hol}(\Omega)$ extends to a full neighborhood of a point $z^1\in L$, then it also extends to a neighborhood of any other point $z^o\in L$.
\et
\bpf
We may assume that $z^o=(0,0)$, $z^1=(|\log\frac14|^{-1},0)$ and that $f$ extends to $B_\delta(z^1)$,  the $\delta$-neighborhood of $z^1$.
Recall that the points $z^o$ and $z^1$ correspond to $\tau=1$ and $\tau=-1$ respectively under the map $\phi_\alpha$. We also remark that 
$$
\phi_\alpha([-\pi,+\pi]\setminus [-e^{-\frac\delta{2\alpha}},e^{\frac\delta{2\alpha}}])\subset B_{\frac\delta2}(z^1).
$$
We deform $h$ by allowing a $\frac\delta2$-bump at $z^1$. Thus, we define
\begin{equation}
\Label{*}
\tilde h=
\begin{cases}
h& \T{ on $[-e^{-\frac\epsilon{2\alpha}},e^{-\frac\epsilon{2\alpha}}]$},
\\
-\frac\delta2&\T{on $[-\pi,\pi]\setminus[-2e^{-\frac\epsilon{2\alpha}},2e^{-\frac\epsilon{2\alpha}}]$},
\end{cases}
\end{equation}
continued smoothly for $e^{-\frac\epsilon{2\alpha}}<|\theta|<2e^{-\frac\epsilon{2\alpha}}$. We 
attach a disc $A_\alpha=(\phi_\alpha,\tilde u_\alpha+i\tilde v_\alpha)$ over $\phi_\alpha$ to the hypersurface defined by $x_2=\tilde h$ according to Proposition~\ref{p1.2}; we have
\begin{equation*}
\begin{split}
\partial_r\tilde u_\alpha&=-\int_{-\pi}^\pi\frac{\tilde u_\alpha}{1-\cos\theta}d\theta
\\
&\geq-2\int_0^{2e^{-\frac\epsilon{2\alpha}}}\frac{h}{1-\cos\theta}d\theta+2\int_{2e^{-\frac\epsilon{2\alpha}}}^\pi\frac\delta{2(1-\cos\theta)}d\theta.
\end{split}
\end{equation*}
Since $s\geq1$, then $\int_0^{2e^{-\frac\epsilon{2\alpha}}}\frac{h}{1-\cos\theta}d\theta\to0$ according to Proposition~\ref{p1.3} (i); thus
$$
\partial_r\tilde u_\alpha>0.
$$
In other terms, $\tilde u_\alpha$ ``points down" at $\tau=1$; in particular,
\begin{equation}
\Label{transversal}
\tilde u_\alpha(1-r)<0\quad\T{for $r<1$ close to $r=1$.}
\end{equation}
We fix $\alpha$ for which \eqref{transversal} is fulfilled and do not keep track of it in the notations which follow. If we replace $\phi$ by $-\epsilon+\phi$, for a fixed $\epsilon$, and substitute in \eqref{*} $-\frac\delta2$ by $-\eta\frac\delta2$ for any $\eta\in[-1,1]$, we get a family of discs
$\{A_\eta\}_\eta=\{(-\epsilon+\phi,\tilde u_\eta+i\tilde v_\eta)\}_\eta$ such that
\begin{equation*}
\begin{cases}
\partial A_\eta\subset \partial\Omega\cup B_\delta(z^1),
\\
\Omega\cup(\underset\eta\cup A_\eta)\T{ contains a neighborhood of $0$}.
\end{cases}
\end{equation*}
Since $f$ extends from the $\partial A_\eta$'s to the $A_\eta$'s by Cauchy's formula, the theorem follows.

\epf

\end{document}